\documentclass[12pt]{article}
\usepackage[top=2.54cm, bottom=2.54cm, left=2.5cm, right=2.5cm]{geometry}
\usepackage[tbtags]{amsmath}
\usepackage{amssymb}
\usepackage{amsthm}
\usepackage{fancyhdr}
\usepackage{latexsym}
\usepackage{mathrsfs}
\usepackage{xcolor}
\usepackage{wasysym}
\usepackage{fancyhdr}
\usepackage{float}
\usepackage{graphicx}
\usepackage[numbers,sort&compress]{natbib}
\usepackage{pict2e,color}
\usepackage{tikz}
\usepackage{enumerate}
\allowdisplaybreaks[4]

\newtheorem{theorem}{\bf Theorem}[section]

\newtheorem{conjecture}[theorem]{Conjecture}

\newtheorem*{claimA1}{Claim A1}
\newtheorem*{claimA2}{Claim A2}
\newtheorem*{claimA3}{Claim A3}
\newtheorem*{claimB1}{Claim B1}
\newtheorem*{claimB2}{Claim B2}
\newtheorem*{claimB3}{Claim B3}

\begin{document}
\title{Towards a conjecture of Birmel\'e-Bondy-Reed on the Erd\H{o}s-P\'osa property of long cycles}
\date{}

\author{Jie Ma\thanks{School of Mathematical Sciences, University of Science and Technology of China, Hefei, Anhui 230026, China.
Email: jiema@ustc.edu.cn. Research supported by the National Key R and D Program of China 2020YFA0713100,
National Natural Science Foundation of China grant 12125106, and Anhui Initiative in Quantum Information Technologies grant AHY150200.}~~~~~
Chunlei Zu\thanks{School of CyberScience and School of Mathematical Sciences, University
of Science and Technology of China, Hefei, Anhui 230026, China.
Email: zucle@mail.ustc.edu.cn}
}

\maketitle

\begin{abstract}
A conjecture of Birmel\'e, Bondy and Reed states that for any integer $\ell\geq 3$,
every graph $G$ without two vertex-disjoint cycles of length at least $\ell$ contains a set of at most $\ell$ vertices which meets all cycles of length at least $\ell$.
They showed the existence of such a set of at most $2\ell+3$ vertices.
This was improved by Meierling, Rautenbach and Sasse to $5\ell/3+29/2$.
Here we present a proof showing that at most $3\ell/2+7/2$ vertices suffice.

\end{abstract}


\section{Introduction}
Let $\mathscr{F}$ be a family of graphs.
For a given graph $G$, a subset $X$ of $V(G)$ is called a \emph{transversal} of $\mathscr{F}$ if the graph $G-X$ contains no member of $\mathscr{F}$.
We say that $\mathscr{F}$ has the \emph{Erd\H{o}s-P\'osa property}, if there is a function $f:\mathbb{N}\rightarrow\mathbb{N}$ such that for every positive integer $k$,
every graph contains either $k$ vertex-disjoint members of $\mathscr{F}$ or a transversal of $\mathscr{F}$ of size at most $f(k)$.
A celebrated result of Erd\H{o}s and P\'osa \cite{EP1965} in 1965 states that the family of all cycles has the Erd\H{o}s-P\'osa property.
Since then it has stimulated a new field of extensive research.

For any integer $\ell\geq3$, let $\mathscr{F}_\ell$ denote the family of cycles of length at least $\ell$.
In 2007, Birmel\'e, Bondy and Reed \cite{BBR2007} first proved that for every $\ell$, $\mathscr{F}_\ell$ has the Erd\H{o}s-P\'osa property.
To be precise, they showed that any graph without $k$ vertex-disjoint cycles in $\mathscr{F}_\ell$ has a transversal of $\mathscr{F}_\ell$ of size at most $O(\ell k^2)$.
The bound of the transversal was improved by Fiorini and Herinckx \cite{FH2014} to $O(\ell k\log k)$.
In 2017, Mousset, Noever, \v{S}kori\'{c} and Weissenberger \cite{MNSW2017} further improved this to $O(\ell k+k\log k)$ and they also provided examples, showing that this is optimal up to the constant factor.

The present paper focuses on the base case $k=2$ of the above problem,
namely, considering graphs without two vertex-disjoint cycles in $\mathscr{F}_\ell$.
As remarked by Birmel\'e, Bondy and Reed \cite{BBR2007}, the case $k=2$ is ``of particular importance''.
Indeed, all proofs of the above papers use inductive arguments.
Birmel\'e, Bondy and Reed \cite{BBR2007} made the following conjecture.

\begin{conjecture}[Birmel\'e, Bondy, and Reed \cite{BBR2007}]\label{Conj}
Let $\ell\geq 3$ and let $G$ be a graph containing no two vertex-disjoint cycles of $\mathscr{F}_\ell$. Then there exists a transversal of $\mathscr{F}_\ell$ of size at most $\ell$.
\end{conjecture}

Note that in view of the complete graph on $2\ell-1$ vertices, the conjectured bound would be best possible.
An early result of Lov\'asz \cite{Lov1965} implied the case $\ell=3$.
Birmel\'e \cite{Bir2003} confirmed the cases $\ell\in\{4,5\}$.
For general $\ell$, Birmel\'e, Bondy and Reed \cite{BBR2007} proved that
there exists a transversal of $\mathscr{F}_\ell$ of size at most $2\ell+3$.
Later, Meierling, Rautenbach and Sasse \cite{MRS2014} improved this to $5\ell/3+29/2$.
Our main result here gives a further improvement as follows.

\begin{theorem}\label{3/2}
Let $\ell\geq3$ be an integer. Let $G$ be a graph containing no two vertex-disjoint cycles of $\mathscr{F}_\ell$. Then there exists a transversal of $\mathscr{F}_\ell$ of size at most $3\ell/2+7/2$.
\end{theorem}

For more references on the Erd\H{o}s-P\'osa property,
we would like to direct interested readers to the survey of Raymond and Thilikos \cite{RT2017}
and \cite{BHJ2019, BJS2018, CJU2020, KKL2020, KK2020, Wei2019} for some recent developments (by no mean of a comprehensive list).
The rest of the paper is organized as follows.
In Section~2 we give the notation, while Section~3 is devoted to the proof of Theorem~\ref{3/2}.

\section{The notation}\label{Notation}
All graphs considered in this paper are finite, undirected and simple.
Let $X$ and $Y$ be subgraphs of a graph $G$.
For a vertex $x$ in $V(G)$, we will use the notation $x\in X$ instead of $x\in V(X)$.
An \emph{$(X,Y)$-path} is a path in $G$ which starts at a vertex of $X$ and ends at a vertex of $Y$ such that no internal vertex is contained in $V(X)\cup V(Y)$.
Here we allow the possibility that $X=Y$.
Let $P$ be a path.
By the {\it length} of $P$, we mean the number of edges in $P$.
If $x$ and $y$ are two vertices of $P$, then $xPy$ denotes the subpath of $P$ with initial vertex $x$ and terminal vertex $y$.
We will reserve the term \emph{disjoint} for \emph{vertex-disjoint}.

Let $C$ be a cycle with a prescribed orientation.
For two vertices $x, y\in V(C)$,
the {\it segment} $xCy$ denotes the unique subpath of $C$ from $x$ to $y$ following the orientation of $C$.
So $xCy$ and $yCx$ are edge-disjoint whose union forms the cycle $C$.
Consider two disjoint $(C,C)$-paths $P$ and $P'$ such that $P$ is between $u$ and $v$ and $P'$ is between $u'$ and $v'$.
We say that $P$ and $P'$ are \emph{parallel} (with respect to $C$) if $u,u',v',v$ appear in the given cyclic order on $C$ and \emph{crossing} (with respect to $C$) otherwise (see Figure 1).

\begin{center}\label{Fig:parallel+crossing}
\begin{tikzpicture}[thin][node distance=1cm,on grid]
\draw(0,0)circle(1.5cm);
\node[circle,inner sep=0.3mm,draw=black,fill=black](x1)at(70:1.5) [label=above:$u$] {};
\node[circle,inner sep=0.3mm,draw=black,fill=black](x2)at(110:1.5) [label=above:$u'$] {};
\node[circle,inner sep=0.3mm,draw=black,fill=black](x3)at(250:1.5) [label=below:$v'$] {};
\node[circle,inner sep=0.3mm,draw=black,fill=black](x4)at(290:1.5) [label=below:$v$] {};
\draw[-](x1)to(x4);
\draw[-](x2)to(x3);
\end{tikzpicture}
\hspace{3cm}
\begin{tikzpicture}[thin][node distance=1cm,on grid]
\draw(0,0)circle(1.5cm);
\node[circle,inner sep=0.3mm,draw=black,fill=black](x1)at(70:1.5) [label=above:$u$] {};
\node[circle,inner sep=0.3mm,draw=black,fill=black](x2)at(110:1.5) [label=above:$u'$] {};
\node[circle,inner sep=0.3mm,draw=black,fill=black](x3)at(250:1.5) [label=below:$v$] {};
\node[circle,inner sep=0.3mm,draw=black,fill=black](x4)at(290:1.5) [label=below:$v'$] {};
\draw[-](x1)to(x3);
\draw[-](x2)to(x4);
\end{tikzpicture}
\end{center}
\begin{center}
\begin{tikzpicture}
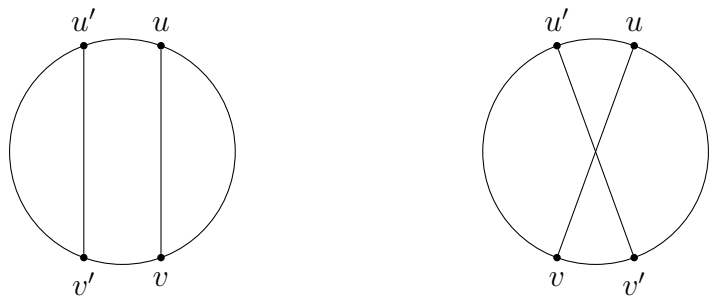
\draw (-3,0) node {Figure 1. parallel and crossing paths};
\end{tikzpicture}
\end{center}

\section{Proof of Theorem \ref{3/2}}\label{Main-Result}
Throughout the rest of this paper, let $\ell\geq 3$ be fixed.
A cycle is called \textbf{long} if it has length at least $\ell$ ( i.e., a cycle in $\mathscr{F}_\ell$) and \textbf{short} otherwise.
We will assume by default that the orientation of any cycle is counterclockwise in all presentations and figures below.

Consider any graph $G$ which contains no two disjoint long cycles.
Our goal is to show that there exists a transversal of $\mathscr{F}_\ell$ of size at most $(3\ell+7)/2$ in $G$.

We now choose two long cycles $C$ and $D$ in $G$ for the coming proof.
Let $C$ be a shortest long cycle of $G$ with length $L$.
It is clear that $C$ intersects every long cycle of $G$, thus $V(C)$ is a transversal of $\mathscr{F}_\ell$.
If $L\leq (3\ell+7)/2$, then the result follows.
So we may assume that $$L>(3\ell+7)/2.$$
We may also assume that there are at least $(3\ell+7)/2\geq 8$ long cycles in $G$ (as otherwise, there is a transversal of $\mathscr{F}_\ell$ of size at most $(3\ell+7)/2$ by taking a vertex from each long cycle).
For every long cycle $D$ of $G$ other than $C$,
let $C_D$ denote a shortest segment of $C$ containing all vertices in $V(C)\cap V(D)$.
Note that $1\leq |V(C_D)|\leq L$.
Choose a long cycle $D$ such that $|V(C_D)|$ is minimum.
With respect to the given orientation of $C$,
we let $x$ and $y$ be the first and last vertices of $C_D$, respectively.
Clearly, $x,y\in V(C)\cap V(D)$.

The rest of the proof will be divided into two cases depending on whether $x=y$ or not.
In each case, using Menger's theorem,
we will find either two disjoint long cycles or a transversal of $\mathscr{F}_\ell$ of size at most $(3\ell+7)/2$, thereby finishing the proof of Theorem \ref{3/2}.

\subsection{The case when $x\neq y$}

Let $X_1$ be the set of $\lceil \ell/2\rceil-1$ vertices of $C$ immediately preceding $x$, and let $X_2$ be the set of $\lceil \ell/2\rceil-1$ vertices of $C$ immediately following $y$. Let $B=C\setminus(X_1\cup X_2\cup V(C_D))$.

We may assume that $G-(X_1\cup X_2\cup\{x,y\})$ contains some long cycle (as otherwise, $X_1\cup X_2\cup\{x,y\}$ is a transversal of $\mathscr{F}_\ell$ of size at most $2\lceil \ell/2\rceil\leq \ell+1$).
Hence every long cycle $D'$ in $G-(X_1\cup X_2\cup\{x,y\})$ intersects $B$ by the minimality of $C_D$.
Let $x_{D'}C_{D'}y_{D'}$ be a shortest segment of $C$ containing $V(B)\cap V(D')$. From now on, choose a long cycle $D'$ such that $|V(x_{D'}Cx)|$ is minimum.

Let $X_3$ be the set of $\lceil \ell/2\rceil-1$ vertices of $C$ immediately preceding $x_{D'}$.
Clearly, $X_1$, $X_2$ and $X_3$ are pairwise disjoint. Otherwise, $X_1\cup X_2\cup X_3\cup\{x,y,x_{D'}\}$ is a transversal of $\mathscr{F}_\ell$.
Since $|X_1\cup X_2\cup X_3\cup\{x,y,x_{D'}\}|\leq\sum_{i=1}^3|X_i|+3=3\lceil \ell/2\rceil\leq3(\ell+1)/2$, we obtain the desired result.
We know that $B\setminus X_3$ consists of two segments of $C$, say $E_1$ and $E_2$.
One is adjacent to $X_1$ and another is adjacent to $X_2$ on $C$. Without loss of generality, we assume that $E_1$ is adjacent to $X_1$ and $E_2$ is adjacent to $X_2$ on $C$.
Note that it is possible that $V(E_1)$ or $V(E_2)$ is empty.

Note that $C\setminus(X_1\cup X_2\cup X_3)$ consists of three segments of $C$, namely $E_1$, $E_2$ and $E_3$ (where $E_3:=C_D$). 
A $(C,C)$-path $P$ with two endpoints $x_0$ and $y_0$ is called a \textbf{special} path between $E_i$ and $E_j$, if $x_0\in V(E_i)$, $y_0\in V(E_j)$ and $i\neq j\in[3]$.

\begin{claimA1}\label{special}
Every special path has length at least $\ell-1$.
\end{claimA1}
\begin{proof}
Let $P$ be a special path between two vertices $x_0$ and $y_0$ of $C$. Let $L_P$ be the length of $P$. Assume by symmetry that $x_0\in V(E_1)$ and $y_0\in V(E_2)$. Since $x_0Cy_0$ has length at least $\ell-1$, $x_0Cy_0\cup y_0Px_0$ forms a long cycle. By the minimality of $C$, $L_P\geq|X_3|+1=\lceil \ell/2\rceil$. Since $y_0Cx_0$ has length at least $\lceil \ell/2\rceil$, we have that $y_0Cx_0\cup x_0Py_0$ is also a long cycle. Thus the length of $P$ is at least the length of $x_0Cy_0$, that is $L_P\geq|X_1|+|X_2|+1=\lceil \ell/2\rceil\times2-1\geq \ell-1$, as desired.
\end{proof}

By the choice of $D'$, we see $D'$ is disjoint from $X_1\cup X_2\cup X_3\cup V(E_2)\cup\{x,y\}$.
Note that $D'$ intersects $D$. It follows that there exists a $(E_1,D\setminus\{x,y\})$-path $sQ_1t$ in $G-(X_1\cup X_2\cup X_3\cup V(E_2)\cup\{x,y\})$, where $s\in V(E_1)$ and $t\in V(D)\setminus\{x,y\}$.

We may assume that there is still a long cycle $D''$ in $G-(X_1\cup X_2\cup X_3\cup\{x,y,x_{D'},t\})$.
This is because that, otherwise, $X_1\cup X_2\cup X_3\cup\{x,y,x_{D'},t\}$ is a transversal of $\mathscr{F}_\ell$ of size at most $\sum_{i=1}^3|X_i|+4=\lceil \ell/2\rceil\times3+1\leq3\ell/2+5/2$.
By the minimality of $C_D$ and the choice of $D'$, we know that $D''$ intersects $E_2$. Moreover, $D''$ intersects $D$.
So there exists a $(E_2,D\setminus\{x,y,t\})$-path $uQ_2v$ in $G-(X_1\cup X_2\cup X_3\cup\{x,y,x_{D'},t\})$, where $u\in V(E_2)$ and $v\in V(D)\setminus\{x,y,t\}$.
We assert that $Q_2\setminus\{u,v\}$ is disjoint from $C\cup sQ_1t$.
Indeed, if not, then there is a special path between $E_1$ and $E_2$ from which it is easy to find a long cycle disjoint from $D$, a contradiction.
Next, we show the following.

\begin{claimA2}\label{v}
$v\in V(tDx)\setminus\{x,t\}$.
\end{claimA2}
\begin{proof}
We have $v\in V(D)\setminus\{x,y,t\}$ and there are three segments of $D\setminus\{x,y,t\}$, namely $xDy\setminus\{x,y\}$, $yDt\setminus\{y,t\}$ and $tDx\setminus\{x,t\}$ (see Figure 2).
Let $C_1:=sCx\cup tDx\cup sQ_1t$. Clearly, $tDx\cup sQ_1t$ contains a special path between $E_1$ and $E_3$.
If $v\in V(xDy)\setminus\{x,y\}$, then $C_2:=yCu\cup uQ_2v\cup vDy$ and $uQ_2v\cup vDy$ contains a special path between $E_2$ and $E_3$,
and if $v\in V(yDt)\setminus\{y,t\}$, then $C_2:=yCu\cup uQ_2v\cup yDv$ and $uQ_2v\cup yDv$ contains a special path between $E_2$ and $E_3$.
By Claim A1, both $C_1$ and $C_2$ are long cycles. So in each case, we find two disjoint long cycles, a contradiction.
\end{proof}

\begin{center}
\begin{tikzpicture}
      \draw [->,thick] (-2.6,2.6) arc (134:158:3cm);
      \draw [draw,thick,color=white] (0,0) ellipse [x radius=2.5cm, y radius=2.5cm];
      \draw [draw,semithick,color=black] (0,-3) ellipse [x radius=1.5cm, y radius=1.5cm];
      \draw [draw,semithick,color=black](1.24,-2.18) arc (-60:240:2.5cm);
      \draw [draw,dashed,thick,color=black](-1.24,-2.18) arc (240:300:2.5cm);
      \draw (-3.5,-2) node {$C$};
      \draw (0,-4.9) node {$D$};
      \draw (-2.37,-1.71) node {$X_1$};
      \draw (2.37,-1.71) node {$X_2$};
      \draw (-0.96,2.7) node {$X_3$};
      \draw (-3.2,0.7) node {$E_1$};
      \draw (2.7,2) node {$E_2$};
      \draw (-1.1,0.95) node {$sQ_1t$};
      \draw (1.15,0.35) node {$uQ_2v$};

   \tikzstyle{every node}=[draw,circle,fill=white,minimum size=3pt,
                            inner sep=0pt]
      \draw [thick](-1.24,-2.18) node (1) [label=below:$x$] {};
      \draw [thick](1.24,-2.18) node (2) [label=below:$y$] {};
      \draw [thick](-1.7,1.82) node (3) [label=145:$x_{D'}$] {};
   \tikzstyle{every node}=[draw,circle,fill=black,minimum size=3pt,
                            inner sep=0pt]
      \draw [thick](-1.48,-2.01) node (4) {};
      \draw [thick](-2.37,-0.8) node (5) {};
      \draw [draw,ultra thick,color=black](-2.37,-0.8) arc (198.5:235:2.5cm);
      \draw [thick](1.48,-2.01) node (6) {};
      \draw [thick](2.37,-0.8) node (7) {};
      \draw [draw,ultra thick,color=black](1.48,-2.01) arc (306.5:340:2.5cm);
      \draw [thick](-1.5,2) node (8) {};
      \draw [thick](0,2.5) node (9) {};
      \draw [draw,ultra thick,color=black](0,2.5) arc (90:127:2.5cm);

      \draw [thick](-2.2,1.2) node (10) [label=left:$s$] {};
      \draw [thick](-0.7,-1.67) node (11) [label=below:$t$] {};
      \draw [draw,thick,color=orange](-0.7,-1.61) arc (10:44.45:5.3cm);
      \draw [thick](1.2,2.2) node (12) [label=above:$u$] {};
      \draw [thick](0.6,-4.37) node (13) [label=below:$v$] {};
      \draw [draw,thick,color=blue](1.17,2.17) arc (150:199.5:7.8cm);
\end{tikzpicture}
\hspace{1.7cm}
\begin{tikzpicture}
      \draw [->,thick] (-2.6,2.6) arc (134:158:3cm);
      \draw [draw,thick,color=white] (0,0) ellipse [x radius=2.5cm, y radius=2.5cm];
      \draw [draw,semithick,color=black] (0,-3) ellipse [x radius=1.5cm, y radius=1.5cm];
      \draw [draw,semithick,color=black](1.24,-2.18) arc (-60:240:2.5cm);
      \draw [draw,dashed,thick,color=black](-1.24,-2.18) arc (240:300:2.5cm);
      \draw (-3.5,-2) node {$C$};
      \draw (0,-4.9) node {$D$};
      \draw (-2.37,-1.71) node {$X_1$};
      \draw (2.37,-1.71) node {$X_2$};
      \draw (-0.96,2.7) node {$X_3$};
      \draw (-3.2,0.7) node {$E_1$};
      \draw (2.7,2) node {$E_2$};
      \draw (-1.1,0.95) node {$sQ_1t$};
      \draw (1.15,0.35) node {$uQ_2v$};

   \tikzstyle{every node}=[draw,circle,fill=white,minimum size=3pt,
                            inner sep=0pt]
      \draw [thick](-1.24,-2.18) node (1) [label=below:$x$] {};
      \draw [thick](1.24,-2.18) node (2) [label=below:$y$] {};
      \draw [thick](-1.7,1.82) node (3) [label=145:$x_{D'}$] {};
   \tikzstyle{every node}=[draw,circle,fill=black,minimum size=3pt,
                            inner sep=0pt]
      \draw [thick](-1.48,-2.01) node (4) {};
      \draw [thick](-2.37,-0.8) node (5) {};
      \draw [draw,ultra thick,color=black](-2.37,-0.8) arc (198.5:235:2.5cm);
      \draw [thick](1.48,-2.01) node (6) {};
      \draw [thick](2.37,-0.8) node (7) {};
      \draw [draw,ultra thick,color=black](1.48,-2.01) arc (306.5:340:2.5cm);
      \draw [thick](-1.5,2) node (8) {};
      \draw [thick](0,2.5) node (9) {};
      \draw [draw,ultra thick,color=black](0,2.5) arc (90:127:2.5cm);

      \draw [thick](-2.2,1.2) node (10) [label=left:$s$] {};
      \draw [thick](-0.7,-1.67) node (11) [label=below:$t$] {};
      \draw [draw,thick,color=orange](-0.7,-1.61) arc (10:44.45:5.3cm);
      \draw [thick](1.2,2.2) node (12) [label=above:$u$] {};
      \draw [thick](0.13,-1.5) node (13) [label=below:$v$] {};
      \draw [draw,thick,color=blue](1.17,2.17) arc (150:177.8:7.8cm);
\end{tikzpicture}
\end{center}
\begin{center}
\begin{tikzpicture}
\draw (-3,0) node {Figure 2. $v\in V(xDy)\setminus\{x,y\}$ and $v\in V(yDt)\setminus\{y,t\}$.};
\end{tikzpicture}
\end{center}

\begin{center}
\begin{tikzpicture}
      \draw [->,thick] (-2.6,2.6) arc (134:158:3cm);
      \draw [draw,thick,color=white] (0,0) ellipse [x radius=2.5cm, y radius=2.5cm];
      \draw [draw,semithick,color=black] (0,-3) ellipse [x radius=1.5cm, y radius=1.5cm];
      \draw [draw,semithick,color=black](1.24,-2.18) arc (-60:240:2.5cm);
      \draw [draw,dashed,thick,color=black](-1.24,-2.18) arc (240:300:2.5cm);
      \draw (-3.5,-2) node {$C$};
      \draw (0,-4.9) node {$D$};
      \draw (-2.37,-1.71) node {$X_1$};
      \draw (2.37,-1.71) node {$X_2$};
      \draw (-0.96,2.7) node {$X_3$};
      \draw (-3.2,0.7) node {$E_1$};
      \draw (2.7,2) node {$E_2$};
      \draw (-1.55,0.2) node {$sQ_1t$};
      \draw (-0.2,1.37) node {$pQq$};
      \draw (1.25,0) node {$p'Q'q'$};

   \tikzstyle{every node}=[draw,circle,fill=white,minimum size=3pt,
                            inner sep=0pt]
      \draw [thick](-1.24,-2.18) node (1) [label=below:$x$] {};
      \draw [thick](1.24,-2.18) node (2) [label=below:$y$] {};
      \draw [thick](-1.7,1.82) node (3) [label=145:$x_{D'}$] {};
   \tikzstyle{every node}=[draw,circle,fill=black,minimum size=3pt,
                            inner sep=0pt]
      \draw [thick](-1.48,-2.01) node (4) {};
      \draw [thick](-2.37,-0.8) node (5) {};
      \draw [draw,ultra thick,color=black](-2.37,-0.8) arc (198.5:235:2.5cm);
      \draw [thick](1.48,-2.01) node (6) {};
      \draw [thick](2.37,-0.8) node (7) {};
      \draw [draw,ultra thick,color=black](1.48,-2.01) arc (306.5:340:2.5cm);
      \draw [thick](-1.5,2) node (8) {};
      \draw [thick](0,2.5) node (9) {};
      \draw [draw,ultra thick,color=black](0,2.5) arc (90:127:2.5cm);

      \draw [thick](-2.2,1.2) node (10) [label=left:$s$] {};
      \draw [thick](0.7,-1.67) node (11) [label=below:$t$] {};
      \draw [draw,thick,color=orange](0.7,-1.61) arc (23.5:67.7:5.3cm);
      \draw [thick](1.2,2.2) node (12) [label=above:$p$] {};
      \draw [thick](-0.7,-1.68) node (13) [label=below:$q$] {};
      \draw [draw,thick,color=purple](1.17,2.17) arc (138:169.5:7.8cm);
      \draw [thick](1.7,1.8) node (12) [label=60:$p'$] {};
      \draw [thick](-0.07,-1.5) node (13) [label=below:$q'$] {};
      \draw [draw,thick,color=blue](1.7,1.8) arc (138:165.5:7.8cm);
\end{tikzpicture}
\hspace{1.7cm}
\begin{tikzpicture}
      \draw [->,thick] (-2.6,2.6) arc (134:158:3cm);
      \draw [draw,thick,color=white] (0,0) ellipse [x radius=2.5cm, y radius=2.5cm];
      \draw [draw,semithick,color=black] (0,-3) ellipse [x radius=1.5cm, y radius=1.5cm];
      \draw [draw,semithick,color=black](1.24,-2.18) arc (-60:240:2.5cm);
      \draw [draw,dashed,thick,color=black](-1.24,-2.18) arc (240:300:2.5cm);
      \draw (-3.5,-2) node {$C$};
      \draw (0,-4.9) node {$D$};
      \draw (-2.37,-1.71) node {$X_1$};
      \draw (2.37,-1.71) node {$X_2$};
      \draw (-0.96,2.7) node {$X_3$};
      \draw (-3.2,0.7) node {$E_1$};
      \draw (2.7,2) node {$E_2$};
      \draw (-1.55,0.2) node {$sQ_1t$};
      \draw (0,1.4) node {$pQq$};
      \draw (1.45,0.7) node {$p'Q'q'$};

   \tikzstyle{every node}=[draw,circle,fill=white,minimum size=3pt,
                            inner sep=0pt]
      \draw [thick](-1.24,-2.18) node (1) [label=below:$x$] {};
      \draw [thick](1.24,-2.18) node (2) [label=below:$y$] {};
      \draw [thick](-1.7,1.82) node (3) [label=145:$x_{D'}$] {};
   \tikzstyle{every node}=[draw,circle,fill=black,minimum size=3pt,
                            inner sep=0pt]
      \draw [thick](-1.48,-2.01) node (4) {};
      \draw [thick](-2.37,-0.8) node (5) {};
      \draw [draw,ultra thick,color=black](-2.37,-0.8) arc (198.5:235:2.5cm);
      \draw [thick](1.48,-2.01) node (6) {};
      \draw [thick](2.37,-0.8) node (7) {};
      \draw [draw,ultra thick,color=black](1.48,-2.01) arc (306.5:340:2.5cm);
      \draw [thick](-1.5,2) node (8) {};
      \draw [thick](0,2.5) node (9) {};
      \draw [draw,ultra thick,color=black](0,2.5) arc (90:127:2.5cm);

      \draw [thick](-2.2,1.2) node (10) [label=left:$s$] {};
      \draw [thick](0.7,-1.67) node (11) [label=below:$t$] {};
      \draw [draw,thick,color=orange](0.7,-1.61) arc (23.5:67.7:5.3cm);
      \draw [thick](1.2,2.2) node (12) [label=above:$p$] {};
      \draw [thick](-0.07,-1.5) node (13) [label=below:$q$] {};
      \draw [draw,thick,color=purple](1.17,2.17) arc (146.5:175:7.8cm);
      \draw [thick](1.7,1.8) node (12) [label=60:$p'$] {};
      \draw [thick](-0.7,-1.68) node (13) [label=below:$q'$] {};
      \draw [draw,thick,color=blue](1.7,1.8) arc (129.5:160.5:7.8cm);
\end{tikzpicture}
\end{center}
\begin{center}
\begin{tikzpicture}
\draw (-3,0) node {Figure 3. Two configurations in the proof of Claim A3.};
\end{tikzpicture}
\end{center}

Now, we see that $uQ_2v$ is a $(E_2,tDx\setminus\{x,t\})$-path in $G-(X_1\cup X_2\cup X_3\cup\{x,y,x_{D'},t\})$ which has no internal vertex in $V(D\cup E_1\cup sQ_1t)$.

\begin{claimA3}\label{1}
One cannot find two disjoint $(E_2,tDx\setminus\{x,t\})$-paths in $G-(X_1\cup X_2\cup X_3\cup\{x,y,x_{D'},t\})$ which has no internal vertex in $V(D\cup E_1\cup sQ_1t)$.
\end{claimA3}
\begin{proof}
Suppose for a contradiction that such two paths exist, say $pQq$ and $p'Q'q'$.
There are two configurations as indicated in Figure 3.
In the left configuration of the figure, 
we have two cycles $C_1:=pCx\cup pQq\cup qDx$ and $C_2:=yC p'\cup p'Q'q'\cup yD q'$.
In the right side, 
we also have two cycles $C_1:=pCs\cup sQ_1t\cup pQq\cup tDq$ and $C_2:=yC p'\cup p'Q'q'\cup q'Dy$.
Using Claim A1, we see that in both cases, $C_1$ and $C_2$ are two disjoint long cycles, a contradiction.
\end{proof}

By Menger's theorem, Claim A3 shows that there is a vertex $z$ meeting all $(E_2,tDx\setminus\{x,t\})$-paths in $G-(X_1\cup X_2\cup X_3\cup\{x,y,x_{D'},t\})$ which has no internal vertex in $V(D\cup E_1\cup sQ_1t)$.
Let $X:=X_1\cup X_2\cup X_3\cup\{x,y,z,x_{D'},t\}$.
Note that $|X|\leq\sum_{i=1}^3|X_i|+5=3\lceil \ell/2\rceil+2\leq (3\ell+7)/2$.
So it suffices to show that $X$ is a transversal of $\mathscr{F}_\ell$.
Suppose not. Then there is a long cycle $D^*$ in $G-X$.
Repeating the same proof as above,
one can find a $(E_2,tDx\setminus\{x,t\})$-path in $G-X$ which has no internal vertex in $V(D\cup E_1\cup sQ_1t)$,
a contradiction to the definition of the vertex $z$.
This completes the proof for the case $x\neq y$.

\subsection{The case when $x=y$}
In this case, clearly we may assume that $G-\{x\}$ contains at least one long cycle.
Every long cycle in $G-\{x\}$ intersects each of the long cycles $C$ and $D$.
Thus there exists at least one $(C,D)$-path in $G-\{x\}$.
We choose a $(C,D)$-path $y'P_0w$ in $G-\{x\}$, where $y'\in V(C)$ and $w\in V(D)$, such that the distance in $C$ between $x$ and $y'$ is minimum.
Without loss of generality, we assume that $xC y'$ is a shortest path in $C$ between $x$ and $y'$.

Let $X_1$ be the set of $\lceil \ell/2\rceil-1$ vertices of $C$ immediately preceding $x$, and let $X_2$ be the set of $\lceil \ell/2\rceil-1$ vertices of $C$ immediately following $y'$. Let $A=xCy'$ and $B=C\setminus(X_1\cup X_2\cup V(xCy'))$.
Since $|X_1\cup X_2\cup\{x,y',w\}|\leq\lceil 2\ell/2\rceil+1\leq \ell+2$,
we may assume that there still is a long cycle $D'$ in $G-(X_1\cup X_2\cup\{x,y',w\})$, which intersects both $C$ and $D$.
If $V(D'\cap C)\subseteq V(xCy')$, then by passing $D'$, one can find a path from $V(xCy')\setminus\{x,y'\}$ to $V(D)\setminus\{x\}$ internally disjoint from $C\cup D$, a contradiction to the definition of $wP_0y'$.
Therefore, every such cycle $D'$ intersects $B$.
Denote $x_{D'}C_{D'}y_{D'}$ to be a shortest segment of $C$ containing $V(B)\cap V(D')$. From now on, choose a long cycle $D'$ such that $|V(x_{D'}Cx)|$ is minimum.
Let $X_3$ be the set of $\lceil \ell/2\rceil-1$ vertices of $C$ immediately preceding $x_{D'}$.
As before, we know that $B\setminus X_3$ consists of two segments of $C$, say $E_1$ and $E_2$.
Without loss of generality, we assume that $E_1$ is adjacent to $X_1$ and $E_2$ is adjacent to $X_2$ on $C$.
Let us call a $(C,C)$-path $P$ with two endpoints $x_0$ and $y_0$ as a \textbf{special} path between $E_i$ and $E_j$, if $x_0\in V(E_i)$, $y_0\in V(E_j)$ and $i\neq j\in[3]$.
We point out that $X_1$, $X_2$ and $X_3$ are pairwise disjoint, and every special path has length at least $\ell-1$.

By the choice of $wP_0y'$, there is no $(A\setminus\{x,y'\},D\setminus\{x,w\})$-path internally disjoint from $C$. It follows from the existence of $D'$ that there is a $(E_1,D\setminus\{x,w\})$-path $sQ_1t$ in $G-(X_1\cup X_2\cup X_3\cup V(E_2)\cup\{x,y',w\})$, where $s\in V(E_1)$ and $t\in V(D)\setminus\{x,w\}$, internally disjoint from $C$, $D$ and $P_0$.

Since $|X_1\cup X_2\cup X_3\cup\{x,y',x_{D'},w\}|\leq 3\lceil \ell/2\rceil+1\leq3\ell/2+5/2$, we may assume that there is a long cycle $D''$ in $G-(X_1\cup X_2\cup X_3\cup\{x,y',x_{D'},w\})$.
By the choice of $y'$ and $x_{D'}$, $D''$ intersects $E_2$ and $D\setminus\{x,w\}$.
So there exists a $(E_2,D\setminus\{x,w\})$-path $uQ_2v$ in $G-(X_1\cup X_2\cup X_3\cup\{x,y',x_{D'},w\})$, where $u\in V(E_2)$ and $v\in V(D)\setminus\{x,w\}$.
We point out that $uQ_2v$ has no internal vertex in $V(C\cup D\cup P_0\cup(Q_1\setminus\{t\}))$
(as otherwise, it is easy to find two disjoint long cycles in $G$ as before; see Figure 4). 
Note that $D\setminus\{x,w\}$ consists of two segments of $D$, i.e., $xDw\setminus\{x,w\}$ and $wDx\setminus\{x,w\}$.

\begin{center}
\begin{tikzpicture}
      \draw [->,thick] (-2.6,2.6) arc (134:158:3cm);
      \draw [draw,semithick,color=black] (0,0) ellipse [x radius=2.5cm, y radius=2.5cm];
      \draw [draw,semithick,color=black] (0,-1.7) ellipse [x radius=0.8cm, y radius=0.8cm];
      \draw (-1.1,-1.5) node {$D$};
      \draw (-3.3,-2) node {$C$};
      \draw (-1.2,-2.7) node {$X_1$};
      \draw (2.9,-0.4) node {$X_2$};
      \draw (-0.96,2.7) node {$X_3$};
      \draw (1.7,-2.4) node {$A$};
      \draw (-3.1,0) node {$E_1$};
      \draw (2,2.3) node {$E_2$};
      \draw (-0.6,0.8) node {$sQ_1t$};
      \draw (1.3,-0.8) node {$wP_0y'$};

   \tikzstyle{every node}=[draw,circle,fill=white,minimum size=3pt,
                            inner sep=0pt]
      \draw [thick](0,-2.5) node (1) [label=below:$x$] {};
      \draw [thick](-1.7,1.82) node (3) [label=145:$x_{D'}$] {};
   \tikzstyle{every node}=[draw,circle,fill=black,minimum size=3pt,
                            inner sep=0pt]
      \draw [thick](-0.28,-2.48) node (4) {};
      \draw [thick](-1.7,-1.83) node (5) {};
      \draw [draw,ultra thick,color=black](-1.7,-1.83) arc (227:262:2.5cm);
      \draw [thick](2.25,-1.1) node (6) {};
      \draw [thick](2.47,0.4) node (7) {};
      \draw [draw,ultra thick,color=black](2.25,-1.1) arc (333.8:368:2.5cm);
      \draw [thick](-1.5,2) node (8) {};
      \draw [thick](0,2.5) node (9) {};
      \draw [draw,ultra thick,color=black](0,2.5) arc (90:127:2.5cm);

      \draw [thick](-2.2,1.2) node (10) [label=left:$s$] {};
      \draw [thick](-0.2,-0.95) node (11) [label=below:$t$] {};
      \draw [draw,thick,color=orange](-0.2,-0.95) arc (27:58.5:5.3cm);
      \draw [thick](2.15,-1.3) node (12) [label=-80:$y'$] {};
      \draw [thick](0.6,-1.2) node (13) [label=235:$w$] {};
      \draw [draw,semithick,color=black](2.15,-1.3) arc (70:102:2.8cm);
\end{tikzpicture}
\end{center}
\begin{center}
\begin{tikzpicture}
\draw (-3,0) node {Figure 4. $wP_0y'$.};
\end{tikzpicture}
\end{center}

\begin{claimB1}\label{tv}
If $t\in V(xDw)\setminus\{x,w\}$, then $v\in V(xDt)\setminus\{x\}$;
if $t\in V(wDx)\setminus\{x,w\}$, then $v\in V(tDx)\setminus\{x\}$.
\end{claimB1}
\begin{proof}
First, consider $t\in V(xDw)\setminus\{x,w\}$ (see Figure 5).
Suppose for a contradiction that $v\notin V(xDt)\setminus\{x\}$.
Then either $v\in V(tDw)\setminus\{t,w\}$ or $v\in V(wDx)\setminus\{x,w\}$.
Let $C_1:=sCx\cup xDt\cup sQ_1t$.
If $v\in V(tDw)\setminus\{t,w\}$, then define $C_2:=y'Cu\cup uQ_2v\cup vDw\cup y'P_0w$; otherwise $v\in V(wDx)\setminus\{x,w\}$, define $C_2:=y'Cu\cup uQ_2v\cup wDv\cup y'P_0w$.
In both cases, $C_1$ and $C_2$ are disjoint long cycles, a contradiction.

\begin{center}
\begin{tikzpicture}
      \draw [->,thick] (-2.6,2.6) arc (134:158:3cm);
      \draw [draw,semithick,color=black] (0,0) ellipse [x radius=2.5cm, y radius=2.5cm];
      \draw [draw,semithick,color=black] (0,-1.7) ellipse [x radius=0.8cm, y radius=0.8cm];
      \draw (-1.1,-1.5) node {$D$};
      \draw (-3.3,-2) node {$C$};
      \draw (-1.2,-2.7) node {$X_1$};
      \draw (2.9,-0.4) node {$X_2$};
      \draw (-0.96,2.7) node {$X_3$};
      \draw (1.7,-2.4) node {$A$};
      \draw (-3.1,0) node {$E_1$};
      \draw (2,2.3) node {$E_2$};
      \draw (-1.2,0.4) node {$sQ_1t$};
      \draw (1.4,0.7) node {$uQ_2v$};
      \draw (1.4,-0.7) node {$P_0$};
      \draw (0,-3.8) node {$v\in V(tDw)\setminus\{t,w\}$};

   \tikzstyle{every node}=[draw,circle,fill=white,minimum size=3pt,
                            inner sep=0pt]
      \draw [thick](0,-2.5) node (1) [label=below:$x$] {};
      \draw [thick](-1.7,1.82) node (3) [label=145:$x_{D'}$] {};
   \tikzstyle{every node}=[draw,circle,fill=black,minimum size=3pt,
                            inner sep=0pt]
      \draw [thick](-0.28,-2.48) node (4) {};
      \draw [thick](-1.7,-1.83) node (5) {};
      \draw [draw,ultra thick,color=black](-1.7,-1.83) arc (227:262:2.5cm);
      \draw [thick](2.25,-1.1) node (6) {};
      \draw [thick](2.47,0.4) node (7) {};
      \draw [draw,ultra thick,color=black](2.25,-1.1) arc (333.8:368:2.5cm);
      \draw [thick](-1.5,2) node (8) {};
      \draw [thick](0,2.5) node (9) {};
      \draw [draw,ultra thick,color=black](0,2.5) arc (90:127:2.5cm);

      \draw [thick](-2.2,1.2) node (10) [label=left:$s$] {};
      \draw [thick](0.8,-1.7) node (11) [label=right:$t$] {};
      \draw [draw,thick,color=orange](0.8,-1.7) arc (0:91:2.9cm);
      \draw [thick](1.2,2.2) node (12) [label=above:$u$] {};
      \draw [thick](0.5,-1.07) node (13) [label=below:$v$] {};
      \draw [draw,thick,color=blue](1.2,2.2) arc (155.5:180:7.8cm);
      \draw [thick](2.15,-1.3) node (12) [label=-80:$y'$] {};
      \draw [thick](0,-0.9) node (13) [label=above:$w$] {};
      \draw [draw,semithick,color=black](2.15,-1.3) arc (56.5:101:2.8cm);
\end{tikzpicture}
\hspace{1.7cm}
\begin{tikzpicture}
      \draw [->,thick] (-2.6,2.6) arc (134:158:3cm);
      \draw [draw,semithick,color=black] (0,0) ellipse [x radius=2.5cm, y radius=2.5cm];
      \draw [draw,semithick,color=black] (0,-1.7) ellipse [x radius=0.8cm, y radius=0.8cm];
      \draw (-1.1,-1.5) node {$D$};
      \draw (-3.3,-2) node {$C$};
      \draw (-1.2,-2.7) node {$X_1$};
      \draw (2.9,-0.4) node {$X_2$};
      \draw (-0.96,2.7) node {$X_3$};
      \draw (1.7,-2.4) node {$A$};
      \draw (-3.1,0) node {$E_1$};
      \draw (2,2.3) node {$E_2$};
      \draw (-1.4,0.55) node {$sQ_1t$};
      \draw (1.15,1) node {$uQ_2v$};
      \draw (1.4,-0.7) node {$P_0$};
      \draw (0,-3.8) node {$v\in V(wDt)\setminus\{t,w\}$};

   \tikzstyle{every node}=[draw,circle,fill=white,minimum size=3pt,
                            inner sep=0pt]
      \draw [thick](0,-2.5) node (1) [label=below:$x$] {};
      \draw [thick](-1.7,1.82) node (3) [label=145:$x_{D'}$] {};
   \tikzstyle{every node}=[draw,circle,fill=black,minimum size=3pt,
                            inner sep=0pt]
      \draw [thick](-0.28,-2.48) node (4) {};
      \draw [thick](-1.7,-1.83) node (5) {};
      \draw [draw,ultra thick,color=black](-1.7,-1.83) arc (227:262:2.5cm);
      \draw [thick](2.25,-1.1) node (6) {};
      \draw [thick](2.47,0.4) node (7) {};
      \draw [draw,ultra thick,color=black](2.25,-1.1) arc (333.8:368:2.5cm);
      \draw [thick](-1.5,2) node (8) {};
      \draw [thick](0,2.5) node (9) {};
      \draw [draw,ultra thick,color=black](0,2.5) arc (90:127:2.5cm);

      \draw [thick](-2.2,1.2) node (10) [label=left:$s$] {};
      \draw [thick](0.8,-1.7) node (11) [label=right:$t$] {};
      \draw [draw,thick,color=orange](0.8,-1.7) arc (0:91:2.9cm);
      \draw [thick](1.2,2.2) node (12) [label=above:$u$] {};
      \draw [thick](-0.5,-1.07) node (13) [label=below:$v$] {};
      \draw [draw,thick,color=blue](1.2,2.2) arc (139:166:7.8cm);
      \draw [thick](2.15,-1.3) node (12) [label=-80:$y'$] {};
      \draw [thick](0,-0.9) node (13) [label=below:$w$] {};
      \draw [draw,semithick,color=black](2.15,-1.3) arc (56.5:101:2.8cm);
\end{tikzpicture}
\end{center}
\begin{center}
\begin{tikzpicture}
\draw (-3,0) node {Figure 5. $t\in V(xDw)\setminus\{x,w\}$.};
\end{tikzpicture}
\end{center}

It remains to consider $t\in V(wDx)\setminus\{x,w\}$ (see Figure 6).
Suppose that $v\notin V(tDx)\setminus\{x\}$.
Then either $v\in V(xDw)\setminus\{x,w\}$ or $v\in V(wDt)\setminus\{w,t\}$.
Let $C_3:=sCx\cup tDx\cup sQ_1t$.
If $v\in V(xDw)\setminus\{x,w\}$, then define $C_4:=y'Cu\cup uQ_2v\cup vDw\cup y'P_0w$;
otherwise $v\in V(wDt)\setminus\{w,t\}$, define $C_4:=y'Cu\cup uQ_2v\cup wDv\cup y'P_0w$.
Again, in both cases, $C_3$ and $C_4$ are two disjoint long cycles, a contradiction.
\end{proof}

\begin{center}
\begin{tikzpicture}
      \draw [->,thick] (-2.6,2.6) arc (134:158:3cm);
      \draw [draw,semithick,color=black] (0,0) ellipse [x radius=2.5cm, y radius=2.5cm];
      \draw [draw,semithick,color=black] (0,-1.7) ellipse [x radius=0.8cm, y radius=0.8cm];
      \draw (-1.1,-1.5) node {$D$};
      \draw (-3.3,-2) node {$C$};
      \draw (-1.2,-2.7) node {$X_1$};
      \draw (2.9,-0.4) node {$X_2$};
      \draw (-0.96,2.7) node {$X_3$};
      \draw (1.7,-2.4) node {$A$};
      \draw (-3.1,0) node {$E_1$};
      \draw (2,2.3) node {$E_2$};
      \draw (-0.8,0.7) node {$sQ_1t$};
      \draw (1.4,0.7) node {$uQ_2v$};
      \draw (1.4,-0.7) node {$P_0$};
      \draw (0,-3.8) node {$v\in V(xDw)\setminus\{x,w\}$};

   \tikzstyle{every node}=[draw,circle,fill=white,minimum size=3pt,
                            inner sep=0pt]
      \draw [thick](0,-2.5) node (1) [label=below:$x$] {};
      \draw [thick](-1.7,1.82) node (3) [label=145:$x_{D'}$] {};
   \tikzstyle{every node}=[draw,circle,fill=black,minimum size=3pt,
                            inner sep=0pt]
      \draw [thick](-0.28,-2.48) node (4) {};
      \draw [thick](-1.7,-1.83) node (5) {};
      \draw [draw,ultra thick,color=black](-1.7,-1.83) arc (227:262:2.5cm);
      \draw [thick](2.25,-1.1) node (6) {};
      \draw [thick](2.47,0.4) node (7) {};
      \draw [draw,ultra thick,color=black](2.25,-1.1) arc (333.8:368:2.5cm);
      \draw [thick](-1.5,2) node (8) {};
      \draw [thick](0,2.5) node (9) {};
      \draw [draw,ultra thick,color=black](0,2.5) arc (90:127:2.5cm);

      \draw [thick](-2.2,1.2) node (10) [label=left:$s$] {};
      \draw [thick](-0.5,-1.06) node (11) [label=below:$t$] {};
      \draw [draw,thick,color=orange](-0.5,-1.06) arc (16:57.5:3.9cm);
      \draw [thick](1.2,2.2) node (12) [label=above:$u$] {};
      \draw [thick](0.5,-1.07) node (13) [label=below:$v$] {};
      \draw [draw,thick,color=blue](1.2,2.2) arc (155.5:180:7.8cm);
      \draw [thick](2.15,-1.3) node (12) [label=-80:$y'$] {};
      \draw [thick](0,-0.9) node (13) [label=above:$w$] {};
      \draw [draw,semithick,color=black](2.15,-1.3) arc (56.5:101:2.8cm);
\end{tikzpicture}
\hspace{1.7cm}
\begin{tikzpicture}
      \draw [->,thick] (-2.6,2.6) arc (134:158:3cm);
      \draw [draw,semithick,color=black] (0,0) ellipse [x radius=2.5cm, y radius=2.5cm];
      \draw [draw,semithick,color=black] (0,-1.7) ellipse [x radius=0.8cm, y radius=0.8cm];
      \draw (-1.1,-1.5) node {$D$};
      \draw (-3.3,-2) node {$C$};
      \draw (-1.2,-2.7) node {$X_1$};
      \draw (2.9,-0.4) node {$X_2$};
      \draw (-0.96,2.7) node {$X_3$};
      \draw (1.7,-2.4) node {$A$};
      \draw (-3.1,0) node {$E_1$};
      \draw (2,2.3) node {$E_2$};
      \draw (-1,0.9) node {$sQ_1t$};
      \draw (1.15,0.8) node {$uQ_2v$};
      \draw (1.4,-0.85) node {$P_0$};
      \draw (0,-3.8) node {$v\in V(wDt)\setminus\{t,w\}$};

   \tikzstyle{every node}=[draw,circle,fill=white,minimum size=3pt,
                            inner sep=0pt]
      \draw [thick](0,-2.5) node (1) [label=below:$x$] {};
      \draw [thick](-1.7,1.82) node (3) [label=145:$x_{D'}$] {};
   \tikzstyle{every node}=[draw,circle,fill=black,minimum size=3pt,
                            inner sep=0pt]
      \draw [thick](-0.28,-2.48) node (4) {};
      \draw [thick](-1.7,-1.83) node (5) {};
      \draw [draw,ultra thick,color=black](-1.7,-1.83) arc (227:262:2.5cm);
      \draw [thick](2.25,-1.1) node (6) {};
      \draw [thick](2.47,0.4) node (7) {};
      \draw [draw,ultra thick,color=black](2.25,-1.1) arc (333.8:368:2.5cm);
      \draw [thick](-1.5,2) node (8) {};
      \draw [thick](0,2.5) node (9) {};
      \draw [draw,ultra thick,color=black](0,2.5) arc (90:127:2.5cm);

\draw [thick](-2.2,1.2) node (10) [label=left:$s$] {};
      \draw [thick](-0.5,-1.06) node (11) [label=below:$t$] {};
      \draw [draw,thick,color=orange](-0.5,-1.06) arc (16:57.5:3.9cm);
      \draw [thick](1.2,2.2) node (12) [label=above:$u$] {};
      \draw [thick](0,-0.9) node (13) [label=below:$v$] {};
      \draw [draw,thick,color=blue](1.2,2.2) arc (146.5:171:7.8cm);
      \draw [thick](2.15,-1.3) node (12) [label=-80:$y'$] {};
      \draw [thick](0.5,-1.1) node (13) [label=below:$w$] {};
      \draw [draw,semithick,color=black](2.15,-1.3) arc (65.5:101:2.8cm);
\end{tikzpicture}
\end{center}
\begin{center}
\begin{tikzpicture}
\draw (-3,0) node {Figure 6. $t\in V(wDx)\setminus\{x,w\}$.};
\end{tikzpicture}
\end{center}

Let $pQq$ and $p'Q'q'$ be disjoint $(E_2,D\setminus\{x,w\})$-paths in $G-(X_1\cup X_2\cup X_3\cup\{x,y',x_{D'},w\})$, where $p,p'\in V(E_2)$ and $q,q'\in V(D)\setminus\{x,w\}$, such that they have no internal vertex in $V(C\cup D\cup P_0\cup(Q_1\setminus\{t\}))$.
Without loss of generality, we assume that $p'$ precedes $p$ on $C$.

\begin{center}
\begin{tikzpicture}
      \draw [->,thick] (-2.6,2.6) arc (134:158:3cm);
      \draw [draw,semithick,color=black] (0,0) ellipse [x radius=2.5cm, y radius=2.5cm];
      \draw [draw,semithick,color=black] (0,-1.7) ellipse [x radius=0.8cm, y radius=0.8cm];
      \draw (-1.1,-1.5) node {$D$};
      \draw (-3.3,-2) node {$C$};
      \draw (-1.2,-2.7) node {$X_1$};
      \draw (2.9,-0.4) node {$X_2$};
      \draw (-0.96,2.7) node {$X_3$};
      \draw (1.7,-2.4) node {$A$};
      \draw (-3.1,0) node {$E_1$};
      \draw (1.9,2.6) node {$E_2$};
      \draw (-0.7,0.7) node {$Q$};
      \draw (0.9,0.5) node {$Q'$};
      \draw (1.4,-0.7) node {$P_0$};
      \draw (0,-3.8) node {$q'\in V(xDw\setminus\{x,w\})$};

   \tikzstyle{every node}=[draw,circle,fill=white,minimum size=3pt,
                            inner sep=0pt]
      \draw [thick](0,-2.5) node (1) [label=below:$x$] {};
      \draw [thick](-1.7,1.82) node (3) [label=145:$x_{D'}$] {};
   \tikzstyle{every node}=[draw,circle,fill=black,minimum size=3pt,
                            inner sep=0pt]
      \draw [thick](-0.28,-2.48) node (4) {};
      \draw [thick](-1.7,-1.83) node (5) {};
      \draw [draw,ultra thick,color=black](-1.7,-1.83) arc (227:262:2.5cm);
      \draw [thick](2.25,-1.1) node (6) {};
      \draw [thick](2.47,0.4) node (7) {};
      \draw [draw,ultra thick,color=black](2.25,-1.1) arc (333.8:368:2.5cm);
      \draw [thick](-1.5,2) node (8) {};
      \draw [thick](0,2.5) node (9) {};
      \draw [draw,ultra thick,color=black](0,2.5) arc (90:127:2.5cm);

      \draw [thick](0.8,2.38) node (10) [label=above:$p$] {};
      \draw [thick](-0.5,-1.06) node (11) [label=-85:$q$] {};
      \draw [draw,thick,color=purple](0.8,2.38) arc (131:186.5:3.9cm);
      \draw [thick](2,1.5) node (12) [label=30:$p'$] {};
      \draw [thick](0.5,-1.07) node (13) [label=-145:$q'$] {};
      \draw [draw,thick,color=blue](2,1.5) arc (139:160.7:7.8cm);
      \draw [thick](2.15,-1.3) node (12) [label=-80:$y'$] {};
      \draw [thick](0,-0.9) node (13) [label=above:$w$] {};
      \draw [draw,semithick,color=black](2.15,-1.3) arc (56.5:101:2.8cm);
\end{tikzpicture}
\hspace{1.7cm}
\begin{tikzpicture}
      \draw [->,thick] (-2.6,2.6) arc (134:158:3cm);
      \draw [draw,semithick,color=black] (0,0) ellipse [x radius=2.5cm, y radius=2.5cm];
      \draw [draw,semithick,color=black] (0,-1.7) ellipse [x radius=0.8cm, y radius=0.8cm];
      \draw (-1.1,-1.5) node {$D$};
      \draw (-3.3,-2) node {$C$};
      \draw (-1.2,-2.7) node {$X_1$};
      \draw (2.9,-0.4) node {$X_2$};
      \draw (-0.96,2.7) node {$X_3$};
      \draw (1.7,-2.4) node {$A$};
      \draw (-3.1,0) node {$E_1$};
      \draw (1.9,2.6) node {$E_2$};
      \draw (-0.7,0.7) node {$Q$};
      \draw (1.3,0.3) node {$Q'$};
      \draw (1.4,-0.85) node {$P_0$};
      \draw (0,-3.8) node {$q'\in V(wDx\setminus\{x,w\})$};

   \tikzstyle{every node}=[draw,circle,fill=white,minimum size=3pt,
                            inner sep=0pt]
      \draw [thick](0,-2.5) node (1) [label=below:$x$] {};
      \draw [thick](-1.7,1.82) node (3) [label=145:$x_{D'}$] {};
   \tikzstyle{every node}=[draw,circle,fill=black,minimum size=3pt,
                            inner sep=0pt]
      \draw [thick](-0.28,-2.48) node (4) {};
      \draw [thick](-1.7,-1.83) node (5) {};
      \draw [draw,ultra thick,color=black](-1.7,-1.83) arc (227:262:2.5cm);
      \draw [thick](2.25,-1.1) node (6) {};
      \draw [thick](2.47,0.4) node (7) {};
      \draw [draw,ultra thick,color=black](2.25,-1.1) arc (333.8:368:2.5cm);
      \draw [thick](-1.5,2) node (8) {};
      \draw [thick](0,2.5) node (9) {};
      \draw [draw,ultra thick,color=black](0,2.5) arc (90:127:2.5cm);

      \draw [thick](0.8,2.38) node (10) [label=above:$p$] {};
      \draw [thick](-0.5,-1.06) node (11) [label=-85:$q$] {};
      \draw [draw,thick,color=purple](0.8,2.38) arc (131:186.5:3.9cm);
      \draw [thick](2,1.5) node (12) [label=45:$p'$] {};
      \draw [thick](0,-0.9) node (13) [label=below:$q'$] {};
      \draw [draw,thick,color=blue](2,1.5) arc (128.5:151.5:7.8cm);
      \draw [thick](2.15,-1.3) node (12) [label=-80:$y'$] {};
      \draw [thick](0.5,-1.1) node (13) [label=below:$w$] {};
      \draw [draw,semithick,color=black](2.15,-1.3) arc (65.5:101:2.8cm);
\end{tikzpicture}
\end{center}
\begin{center}
\begin{tikzpicture}
\draw (-3,0) node {Figure 7. $pQq$ and $p'Q'q'$ are parallel.};
\end{tikzpicture}
\end{center}

\begin{claimB2}\label{q}
If $pQq$ and $p'Q'q'$ are parallel, then $q,q'\in V(xDw)\setminus\{x,w\}$;
if $pQq$ and $p'Q'q'$ are crossing, then $q,q'\in V(wDx)\setminus\{x,w\}$;
\end{claimB2}
\begin{proof}
If $pQq$ and $p'Q'q'$ are parallel (see Figure 7), then $q'$ precedes $q$ on $D$.
Suppose for a contradiction that at least one of $q$ and $q'$ is not in $V(xDw)\setminus\{x,w\}$.
Then $q\in V(wDx)\setminus\{x,w\}$.
Let $C_1:=pCx\cup qDx\cup pQq$. If $q'\in V(xDw)\setminus\{x,w\}$, define $C_2:=y'Cp'\cup p'Q'q'\cup q'Dw\cup y'P_0w$;
if $q'\in V(wDx)\setminus\{x,w\}$, define $C_2:=y'Cp'\cup p'Q'q'\cup wDq'\cup y'P_0w$.
In both cases, $C_1$ and $C_2$ are disjoint long cycles, a contradiction.

\begin{center}
\begin{tikzpicture}
      \draw [->,thick] (-2.6,2.6) arc (134:158:3cm);
      \draw [draw,semithick,color=black] (0,0) ellipse [x radius=2.5cm, y radius=2.5cm];
      \draw [draw,semithick,color=black] (0,-1.7) ellipse [x radius=0.8cm, y radius=0.8cm];
      \draw (-1.1,-1.5) node {$D$};
      \draw (-3.3,-2) node {$C$};
      \draw (-1.2,-2.7) node {$X_1$};
      \draw (2.9,-0.4) node {$X_2$};
      \draw (-0.96,2.7) node {$X_3$};
      \draw (1.7,-2.4) node {$A$};
      \draw (-3.1,0) node {$E_1$};
      \draw (1.9,2.6) node {$E_2$};
      \draw (1.1,0) node {$Q$};
      \draw (-0.4,0.3) node {$Q'$};
      \draw (1.4,-0.7) node {$P_0$};
      \draw (0,-3.8) node {$q'\in V(wDx)\setminus\{x,w\}$};

   \tikzstyle{every node}=[draw,circle,fill=white,minimum size=3pt,
                            inner sep=0pt]
      \draw [thick](0,-2.5) node (1) [label=below:$x$] {};
      \draw [thick](-1.7,1.82) node (3) [label=145:$x_{D'}$] {};
   \tikzstyle{every node}=[draw,circle,fill=black,minimum size=3pt,
                            inner sep=0pt]
      \draw [thick](-0.28,-2.48) node (4) {};
      \draw [thick](-1.7,-1.83) node (5) {};
      \draw [draw,ultra thick,color=black](-1.7,-1.83) arc (227:262:2.5cm);
      \draw [thick](2.25,-1.1) node (6) {};
      \draw [thick](2.47,0.4) node (7) {};
      \draw [draw,ultra thick,color=black](2.25,-1.1) arc (333.8:368:2.5cm);
      \draw [thick](-1.5,2) node (8) {};
      \draw [thick](0,2.5) node (9) {};
      \draw [draw,ultra thick,color=black](0,2.5) arc (90:127:2.5cm);

      \draw [thick](0.8,2.38) node (10) [label=above:$p$] {};
      \draw [thick](0.5,-1.07) node (13) [label=-145:$q$] {};
      \draw [draw,thick,color=purple](0.8,2.38) arc (12:-21.5:5.9cm);
      \draw [thick](2,1.5) node (12) [label=30:$p'$] {};
      \draw [thick](-0.5,-1.06) node (11) [label=-85:$q'$] {};
      \draw [draw,thick,color=blue](2,1.5) arc (96:174.5:2.8cm);
      \draw [thick](2.15,-1.3) node (12) [label=-80:$y'$] {};
      \draw [thick](0,-0.9) node (13) [label=above:$w$] {};
      \draw [draw,semithick,color=black](2.15,-1.3) arc (56.5:101:2.8cm);
\end{tikzpicture}
\hspace{1.7cm}
\begin{tikzpicture}
      \draw [->,thick] (-2.6,2.6) arc (134:158:3cm);
      \draw [draw,semithick,color=black] (0,0) ellipse [x radius=2.5cm, y radius=2.5cm];
      \draw [draw,semithick,color=black] (0,-1.7) ellipse [x radius=0.8cm, y radius=0.8cm];
      \draw (-1.1,-1.5) node {$D$};
      \draw (-3.3,-2) node {$C$};
      \draw (-1.2,-2.7) node {$X_1$};
      \draw (2.9,-0.4) node {$X_2$};
      \draw (-0.96,2.7) node {$X_3$};
      \draw (1.7,-2.4) node {$A$};
      \draw (-3.1,0) node {$E_1$};
      \draw (1.9,2.6) node {$E_2$};
      \draw (0.7,1.5) node {$Q$};
      \draw (1.8,0.85) node {$Q'$};
      \draw (1.7,-0.8) node {$P_0$};
      \draw (0,-3.8) node {$q'\in V(xDw)\setminus\{x,w\}$};

   \tikzstyle{every node}=[draw,circle,fill=white,minimum size=3pt,
                            inner sep=0pt]
      \draw [thick](0,-2.5) node (1) [label=below:$x$] {};
      \draw [thick](-1.7,1.82) node (3) [label=145:$x_{D'}$] {};
   \tikzstyle{every node}=[draw,circle,fill=black,minimum size=3pt,
                            inner sep=0pt]
      \draw [thick](-0.28,-2.48) node (4) {};
      \draw [thick](-1.7,-1.83) node (5) {};
      \draw [draw,ultra thick,color=black](-1.7,-1.83) arc (227:262:2.5cm);
      \draw [thick](2.25,-1.1) node (6) {};
      \draw [thick](2.47,0.4) node (7) {};
      \draw [draw,ultra thick,color=black](2.25,-1.1) arc (333.8:368:2.5cm);
      \draw [thick](-1.5,2) node (8) {};
      \draw [thick](0,2.5) node (9) {};
      \draw [draw,ultra thick,color=black](0,2.5) arc (90:127:2.5cm);

      \draw [thick](0.8,2.38) node (10) [label=above:$p$] {};
      \draw [thick](0.8,-1.8) node (11) [label=-85:$q$] {};
      \draw [draw,thick,color=purple](0.8,2.38) arc (20.8:-20:5.9cm);
      \draw [thick](2,1.5) node (12) [label=30:$p'$] {};
      \draw [thick](0.5,-1.07) node (13) [label=-145:$q'$] {};
      \draw [draw,thick,color=blue](2,1.5) arc (126.5:172:3.8cm);
      \draw [thick](2.15,-1.3) node (12) [label=-80:$y'$] {};
      \draw [thick](0,-0.9) node (13) [label=above:$w$] {};
      \draw [draw,semithick,color=black](2.15,-1.3) arc (56.5:101:2.8cm);
\end{tikzpicture}
\end{center}
\begin{center}
\begin{tikzpicture}
\draw (-3,0) node {Figure 8. $pQq$ and $p'Q'q'$ are crossing.};
\end{tikzpicture}
\end{center}

If $pQq$ and $p'Q'q'$ are crossing (see Figure 8), then $q$ precedes $q'$ on $D$.
Suppose for a contradiction that at least one of $q$ and $q'$ is not in $V(wDx)\setminus\{x,w\}$.
Then $q\in V(xDw)\setminus\{x,w\}$.
Let $C_3:=pCx\cup xDq\cup pQq$. If $q'\in V(wDx)\setminus\{x,w\}$, let $C_4:=y'Cp'\cup p'Q'q'\cup wDq'\cup y'P_0w$;
otherwise $q'\in V(xDw)\setminus\{x,w\}$, let $C_4:=y'Cp'\cup p'Q'q'\cup q'Dw\cup y'P_0w$.
In both cases, $C_3$ and $C_4$ are disjoint long cycles, a contradiction.
\end{proof}

\begin{claimB3}\label{11}
One cannot find two disjoint $(E_2,D\setminus\{x,w\})$-paths in $G-(X_1\cup X_2\cup X_3\cup\{x,y',x_{D'},w\})$ which has no internal vertex in $V(C\cup D\cup P_0\cup(Q_1\setminus\{t\}))$.
\end{claimB3}
\begin{proof}
Suppose for a contradiction that such paths exist, say $pQq$ and $p'Q'q'$, where $p,p'\in V(E_2)$, $q,q'\in V(D)\setminus\{x,w\}$ and $p'$ precedes $p$ on $C$.
By Claims B1 and B2, there are two configurations (see Figure 9). Be aware that there might be $q=t$.
In the left configuration of Figure 9, $pQq$, $p'Q'q'$ and $sQ_1t$ are pairwise crossing, and let $C_1:=pCs\cup sQ_1t\cup tDq\cup pQq$ and $C_2:=xCp'\cup p'Q'q'\cup q'Dx$.
In the right configuration, $pQq$, $p'Q'q'$ and $sQ_1t$ are pairwise parallel, and let $C_1:=pCs\cup sQ_1t\cup qDt\cup pQq$ and $C_2:=xCp'\cup p'Q'q'\cup xDq'$.
It is easy to check that $C_1$ and $C_2$ are disjoint long cycles in both cases.
\end{proof}

\begin{center}
\begin{tikzpicture}
      \draw [->,thick] (-2.6,2.6) arc (134:158:3cm);
      \draw [draw,semithick,color=black] (0,0) ellipse [x radius=2.5cm, y radius=2.5cm];
      \draw [draw,semithick,color=black] (0,-1.7) ellipse [x radius=0.8cm, y radius=0.8cm];
      \draw (-1.1,-1.5) node {$D$};
      \draw (-3.3,-2) node {$C$};
      \draw (-1.2,-2.7) node {$X_1$};
      \draw (2.9,-0.4) node {$X_2$};
      \draw (-0.96,2.7) node {$X_3$};
      \draw (1.7,-2.4) node {$A$};
      \draw (-3.1,0) node {$E_1$};
      \draw (1.9,2.6) node {$E_2$};
      \draw (-1,0.75) node {$Q_1$};
      \draw (0.5,1.7) node {$Q$};
      \draw (1.5,1) node {$Q'$};
      \draw (1.5,-1.1) node {$P_0$};

   \tikzstyle{every node}=[draw,circle,fill=white,minimum size=3pt,
                            inner sep=0pt]
      \draw [thick](0,-2.5) node (1) [label=below:$x$] {};
      \draw [thick](-1.7,1.82) node (3) [label=145:$x_{D'}$] {};
   \tikzstyle{every node}=[draw,circle,fill=black,minimum size=3pt,
                            inner sep=0pt]
      \draw [thick](-0.28,-2.48) node (4) {};
      \draw [thick](-1.7,-1.83) node (5) {};
      \draw [draw,ultra thick,color=black](-1.7,-1.83) arc (227:262:2.5cm);
      \draw [thick](2.25,-1.1) node (6) {};
      \draw [thick](2.47,0.4) node (7) {};
      \draw [draw,ultra thick,color=black](2.25,-1.1) arc (333.8:368:2.5cm);
      \draw [thick](-1.5,2) node (8) {};
      \draw [thick](0,2.5) node (9) {};
      \draw [draw,ultra thick,color=black](0,2.5) arc (90:127:2.5cm);

      \draw [thick](-2.2,1.2) node (10) [label=left:$s$] {};
      \draw [thick](0.65,-1.2) node (11) [label=right:$t$] {};
      \draw [draw,thick,color=orange](0.65,-1.2) arc (2:97:2.5cm);
      \draw [thick](0.8,2.38) node (10) [label=above:$p$] {};
      \draw [thick](0.2,-0.93) node (13) [label=-145:$q$] {};
      \draw [draw,thick,color=purple](0.8,2.38) arc (6.5:-26.2:5.9cm);
      \draw [thick](2,1.5) node (12) [label=30:$p'$] {};
      \draw [thick](-0.5,-1.06) node (11) [label=-85:$q'$] {};
      \draw [draw,thick,color=blue](2,1.5) arc (96:174.5:2.8cm);
      \draw [thick](2.15,-1.3) node (12) [label=-80:$y'$] {};
      \draw [thick](0.8,-1.7) node (13) [label=left:$w$] {};
      \draw [draw,semithick,color=black](2.15,-1.3) arc (91.5:120:2.8cm);
\end{tikzpicture}
\hspace{1.7cm}
\begin{tikzpicture}
      \draw [->,thick] (-2.6,2.6) arc (134:158:3cm);
      \draw [draw,semithick,color=black] (0,0) ellipse [x radius=2.5cm, y radius=2.5cm];
      \draw [draw,semithick,color=black] (0,-1.7) ellipse [x radius=0.8cm, y radius=0.8cm];
      \draw (-1.1,-1.5) node {$D$};
      \draw (-3.3,-2) node {$C$};
      \draw (-1.2,-2.7) node {$X_1$};
      \draw (2.9,-0.4) node {$X_2$};
      \draw (-0.96,2.7) node {$X_3$};
      \draw (1.7,-2.4) node {$A$};
      \draw (-3.1,0) node {$E_1$};
      \draw (1.9,2.6) node {$E_2$};
      \draw (-1.2,0.5) node {$Q_1$};
      \draw (0.7,1.1) node {$Q$};
      \draw (1.6,0.7) node {$Q'$};
      \draw (1.7,-1) node {$P_0$};

   \tikzstyle{every node}=[draw,circle,fill=white,minimum size=3pt,
                            inner sep=0pt]
      \draw [thick](0,-2.5) node (1) [label=below:$x$] {};
      \draw [thick](-1.7,1.82) node (3) [label=145:$x_{D'}$] {};
   \tikzstyle{every node}=[draw,circle,fill=black,minimum size=3pt,
                            inner sep=0pt]
      \draw [thick](-0.28,-2.48) node (4) {};
      \draw [thick](-1.7,-1.83) node (5) {};
      \draw [draw,ultra thick,color=black](-1.7,-1.83) arc (227:262:2.5cm);
      \draw [thick](2.25,-1.1) node (6) {};
      \draw [thick](2.47,0.4) node (7) {};
      \draw [draw,ultra thick,color=black](2.25,-1.1) arc (333.8:368:2.5cm);
      \draw [thick](-1.5,2) node (8) {};
      \draw [thick](0,2.5) node (9) {};
      \draw [draw,ultra thick,color=black](0,2.5) arc (90:127:2.5cm);

      \draw [thick](-2.2,1.2) node (10) [label=left:$s$] {};
      \draw [thick](0,-0.9) node (11) [label=below:$t$] {};
      \draw [draw,thick,color=orange](0,-0.9) arc (9:83:2.5cm);
      \draw [thick](0.8,2.38) node (10) [label=above:$p$] {};
      \draw [thick](0.5,-1.07) node (11) [label=-145:$q$] {};
      \draw [draw,thick,color=purple](0.8,2.38) arc (12.5:-21.5:5.9cm);
      \draw [thick](2,1.5) node (12) [label=30:$p'$] {};
      \draw [thick](0.8,-1.8) node (13) [label=-85:$q'$] {};
      \draw [draw,thick,color=blue](2,1.5) arc (-2.5:-37:5.9cm);
      \draw [thick](2.15,-1.3) node (12) [label=-80:$y'$] {};
      \draw [thick](-0.5,-1.05) node (13) [label=below:$w$] {};
      \draw [draw,semithick,color=black](2.15,-1.3) arc (22:150:1.5cm);
\end{tikzpicture}
\end{center}
\begin{center}
\begin{tikzpicture}
\draw (-3,0) node {Figure 9. Two configurations in the proof of Claim B3.};
\end{tikzpicture}
\end{center}

By Menger's theorem, Claim B3 shows that there is a vertex $z$ meeting all $(E_2,D\setminus\{x,w\})$-paths in $G-(X_1\cup X_2\cup X_3\cup\{x,y',x_{D'},w\})$ which has no internal vertex in $V(C\cup D\cup P_0\cup(Q_1\setminus\{t\}))$.
Let $X:=X_1\cup X_2\cup X_3\cup\{x,y',z,x_{D'},w\}$. Note that $|X|\leq\sum_{i=1}^3|X_i|+5=\lceil \ell/2\rceil\times3+2\leq3\ell/2+7/2$.
So it is enough to show that $X$ is a transversal of $\mathscr{F}_\ell$.
Suppose not, then there is a long cycle $D^*$ in $G-X$.
As the same proof, one can show that there exists a $(E_2,D\setminus\{x,w\})$-path in $G-X$ which has no internal vertex in $V(C\cup D\cup P_0\cup(Q_1\setminus\{t\}))$.
This is a contradiction to the definition of the vertex $z$.
We have completed the proof of the case $x=y$ and thereby the proof of Theorem \ref{3/2}.\qed


\begin{thebibliography}{99}
\bibitem{Bir2003} E. Birmel\'{e}, Th\`{e}se de doctorat, Universit\'{e} de Lyon 1, 2003.

\bibitem{BBR2007} E. Birmel\'{e}, J. A. Bondy, and B. A. Reed, The Erd\H{o}s-P\'{o}sa property for long circuits, Combinatorica 27 (2) (2007), 135--145.

\bibitem{BHJ2019} H. Bruhn, M. Heinlein, and F. Joos, Long cycles have the edge-Erd\H{o}s-P\'{o}sa property, Combinatorica 39 (2019), no.1, 1--36.

\bibitem{BJS2018} H. Bruhn, F. Joos, and O. Schaudt, Long cycles through prescribed vertices have the Erd\H{o}s-P\'{o}sa property, J. Graph Theory 87 (2018), no.3, 275--284.

\bibitem{CJU2020} W. Cames Van Batenburg, G. Joret, and A. Ulmer, Erd\H{o}s-P\'{o}sa from ball packing, SIAM J. Discrete Math. 34 (2020), no.3, 1609--1619.

\bibitem{EP1965} P. Erd\H{o}s and L. P\'{o}sa, On independent circuits contained in a graph, Canad. J. Math. 17 (1965), 347--352.

\bibitem{FH2014} S. Fiorini and A. Herinckx, A tighter Erd\H{o}s-P\'{o}sa function for long cycles, J. Graph Theory 77 (2014), no.2, 111--116.

\bibitem{KKL2020} M. Kang, O. Kwon, and M. Lee, Graphs without two vertex-disjoint $S$-cycles. Discrete Math. 343 (2020), no.10, 111997, 18 pp.

\bibitem{KK2020} E. Kim and O. Kwon, Erd\H{o}s-P\'{o}sa property of chordless cycles and its applications. J. Combin. Theory Ser. B 145 (2020), 65--112.

\bibitem{Lov1965} L. Lov\'{a}sz, On graphs not containing independent circuits (Hungarian), Mat. Lapok 16 (1965), 289--299.

\bibitem{MRS2014} D. Meierling, D. Rautenbach, and T. Sasse, The Erd\H{o}s-P\'{o}sa property for long circuits, J. Graph Theory 77 (2014), no.4, 251--259.

\bibitem{MNSW2017} F. Mousset, A. Noever, N. \v{S}kori\'{c}, and F. Weissenberger, A tight Erd\H{o}s-P\'{o}sa function for long cycles, J. Combin. Theory Ser. B 125 (2017), 21--32.

\bibitem{RT2017} J.-F. Raymond and D. M. Thilikos, Recent techniques and results on the Erd\H{o}s-P\'{o}sa property, Discrete Appl. Math. 231 (2017), 25--43.

\bibitem{Wei2019} D. Wei{\ss}auer, In absence of long chordless cycles, large tree-width becomes a local phenomenon. J. Combin. Theory Ser. B 139 (2019), 342--352.

\end{thebibliography}
\end{document}